\documentclass[a4paper, 11pt, twoside]{article}
\usepackage{mfpaperstuff2}
\begin{document}
\newcommand\shift[1]{\vec{#1}}
\newcommand\bicl{$2$-clique\xspace}
\newcommand\bicls{$2$-cliques\xspace}
\newcommand\seqp[1]{\calp_{#1}}
\newcommand\prmp[1]{\calr_{#1}}
\newcommand\intn[1]{\calq(#1)}
\newenvironment{thmcx}[1]{\begin{thmciting}{\textup{\textbf{\cite{#1}.}}\ }}{\end{thmciting}}
\medmuskip=3mu plus 1mu minus 1mu
\newcommand\bc{$2$-chain\xspace}
\newcommand\bcs{$2$-chains\xspace}
\newcommand\cpt[2]{#1^{(#2)}}
\newcommand\nd[2]{\draw(#1)node[minimum size=2.5mm,inner sep=0pt,draw,fill,circle](#2){};}
\newcommand\ndi[2]{\draw(#1)node[inner,minimum size=2.5mm,inner sep=0pt,draw,fill,circle](#2){};}
\newcommand\cls\trianglelefteqslant
\newcommand\cl\vartriangleleft
\newcommand\clsp{\cls^+}
\newcommand\ncls{\not\trianglelefteqslant}
\newcommand\pmi{P^-}
\newcommand\rev[1]{\operatorname{rev}#1}
\newcommand\dual{^\circ}
\newcommand\ins{opposed\xspace}
\newcommand\insness{opposedness\xspace}
\newcommand\splice\boxbslash
\newcommand\im{\operatorname{im}}
\newcommand\nmo{n{-}1}
\newcommand\nmt{n{-}2}
\newcommand\ipo{i{+}1}
\newcommand\ul\underline
\Yvcentermath1

\title{$2$-chains: an interesting family of posets}

\toptitle

\begin{abstract}
We introduce a new family of finite posets which we call \emph{\bcs}. These first arose in the study of $0$-Hecke algebras, but they admit a variety of different characterisations. We give these characterisations, prove that they are equivalent and derive some numerical results concerning \bcs.
\end{abstract}


\section{Introduction}\label{introsec}

In this paper we introduce a family of finite partially ordered sets which we call \emph{\bcs}. These are remarkable in that they can be described in several quite different ways. Some of these descriptions are by direct construction, and some by restrictive properties. \bcs originally arose in the author's paper \cite{0hecke} on $0$-Hecke algebras, where they appear as structure posets for certain induced modules for Hecke algebras of type A. Apart from this, however, \bcs appear not to have been considered before. The fact that they admit so many different descriptions suggests that they should be studied further, and we hope that this paper is the start of an interesting new avenue in poset theory.

The author must express his gratitude to Jeremy Rickard. Early in this work, the author posted a question about these posets on \textit{MathOverflow}, and Jeremy's answer \cite{mo} provided some key results which underpin a lot of this paper, and without which most of this work would not have been done. We indicate later exactly which parts are due to Jeremy.

We give a brief outline of the structure of the paper. We begin in \cref{bcsec} with a definition of \bcs. In \cref{maxlsec} we consider maximal elements, giving a reduction theorem which allows all \bcs to be constructed recursively. In \cref{consec} we give a direct construction of posets from binary sequences, and show that these give all \bcs. In \cref{splicesec} we introduce a method of joining posets together which we call \emph{splicing}, and show that this can be used to join smaller \bcs to make larger ones. In this way, we give another construction of all \bcs as splices of certain simple \bcs. In \cref{graphsec} we consider a connection with graphs: writing down a definition for graphs analogous to the definition of \bcs leads to another characterisation of \bcs, in terms of their incomparability graphs. In \cref{coxsec} we return to the connection with Coxeter groups of type A, showing how to characterise \bcs in terms of Coxeter elements of symmetric groups. In \cref{covsec} we show how to count covers and linear extensions in \bcs. In \cref{gensec,infinitesec} we consider some possible generalisations.

\section{\bcs}\label{bcsec}

In this section we will introduce \bcs. We begin by recalling some basic definitions for posets.

As usual, a \emph{partial order} is a reflexive asymmetric transitive binary relation. A \emph{partially ordered set} or \emph{poset} is a set $P$ with a given partial order $\cls$. We write $p\cl q$ to mean $p\cls q$ and $p\neq q$. Given $p,q\in P$ with $p\cls q$, we say that $q$ is \emph{above} $p$, and that $p$ is \emph{below} $q$. We say that $q$ \emph{covers} $p$ if $p\cl q$ and there is no $r\in P$ with $p\cl r\cl q$. We say that $p,q$ are \emph{incomparable} if $p\ncls q\ncls p$. The order $\cls$ is a \emph{total order} if for all $p,q\in P$ either $p\cls q$ or $q\cls p$. If $Q\subseteq P$, then $\cls$ induces a partial order on $Q$, which will also be denoted $\cls$. $Q$ is called:
\begin{itemize}
\item
a \emph{chain} if this induced partial order is a total order;
\item
an \emph{antichain} if $p\ncls q$ for all $p,q\in Q$;
\item
an \emph{ideal} (or \emph{down-set}) in $P$ if $p\ncls q$ whenever $q\in Q$ and $p\in P\setminus Q$.
\end{itemize}
The \emph{dual} poset $P^\circ$ has the same underlying set, with $p\cls q$ in $P^\circ$ \iff $q\cls p$ in $P$.

A \emph{refinement} of $\cls$ is a partial order $\clsp$ such that $p\clsp q$ whenever $p\cls q$; if there is at least one pair $p,q$ with $p\clsp q$ but $p\ncls q$, then $\clsp$ is a \emph{proper} refinement of $\cls$.

An element $p\in P$ is \emph{maximal} if there is no $q\in P$ with $p\cl q$. $p$ is the \emph{greatest} element of $P$ if $q\cls p$ for all $q\in Q$. \emph{Minimal} and \emph{least} elements are defined similarly.

In this paper, a partial order will almost always be written as $\cls$, and if we say that $P$ is a poset without specifying the partial order, we mean $(P,\cls)$ is a poset. We will use terms such as \emph{chain}, \emph{maximal} and \emph{above} without explicit reference to a particular partial order, and it should be understood that we are referring to $\cls$. If we use such terms with respect to any other partial order, we will be explicit about which order we are referring to.

An \emph{isomorphism} between two posets $P$ and $R$ is a bijection $\phi:P\to R$ such that for $p,q\in P$ we have $p\cls q$ \iff $\phi(p)\cls\phi(q)$. $P$ and $R$ are \emph{isomorphic} (written $P\cong R$) if there is at least one isomorphism from $P$ to $R$. An \emph{automorphism} of $P$ is an isomorphism from $P$ to $P$.

The \emph{Hasse diagram} of a poset $P$ is a graph drawn in the plane with vertices corresponding to the elements of $P$, with an edge drawn from $p$ to $q$ with strictly increasing $y$-coordinate whenever $q$ covers $p$.

\medskip
For most of this paper, we only consider \emph{finite} posets. We make some comments on the infinite case in \cref{infinitesec}; until then, we take ``poset'' to mean ``finite poset''.

\smallskip
Now we can give our main definition. Suppose $P$ is a poset. We say that $P$ is a \emph{\bc} if the following two conditions are satisfied.
\begin{enumerate}[ref=(\theenumi)]
\item\label{bc1}
There is a unique way to write $P$ as the union of two chains.
\item
$\cls$ is maximal subject to \ref{bc1}, i.e.\ if $\clsp$ is a proper refinement of $\cls$, then there is more than one way to write $P$ as the union of two $\clsp$-chains.
\end{enumerate}

Note that when we say ``union'' in \ref{bc1}, we simply mean set union; so it is permissible to have an element of one chain lying above an element of the other.

There are four isomorphism classes of $5$-element \bcs, given by the following Hasse diagrams.

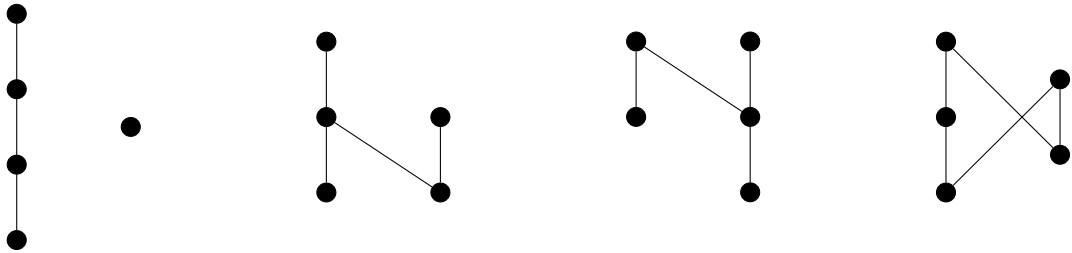
\[
\begin{tikzpicture}[xscale=-1,baseline=.5cm]
\nd{0,1.5}a
\nd{1.5,0}b
\nd{1.5,1}c
\nd{1.5,2}d
\nd{1.5,3}e
\draw(b)--(e);
\end{tikzpicture}
\qquad\qquad\qquad
\begin{tikzpicture}[xscale=-1]
\nd{0,0}a
\nd{0,1}b
\nd{1.5,0}c
\nd{1.5,1}d
\nd{1.5,2}e
\draw(b)--(a)--(d);
\draw(c)--(e);
\end{tikzpicture}
\qquad\qquad\qquad
\begin{tikzpicture}[yscale=-1]
\nd{0,0}a
\nd{0,1}b
\nd{1.5,0}c
\nd{1.5,1}d
\nd{1.5,2}e
\draw(b)--(a)--(d);
\draw(c)--(e);
\end{tikzpicture}
\qquad\qquad\qquad
\begin{tikzpicture}[xscale=-1]
\nd{0,.5}a
\nd{0,1.5}b
\nd{1.5,0}c
\nd{1.5,1}d
\nd{1.5,2}e
\draw(b)--(a)--(e)--(c)--(b);
\end{tikzpicture}
\]
The first example can be generalised to give a \bc of arbitrary size: just take a poset consisting of an arbitrary chain $C$, together with an element $p$ not lying above or below any element of $C$.

For another class of \bcs of arbitrary size, take $n$ to be a positive integer, set $\intn n=\{1,\dots,n\}$, and define a partial order $\cls$ on $\intn n$ by putting $i\cl j$ \iff $i\ls j-2$ in the usual order on $\{1,\dots,n\}$. For example, $\intn8$ has the following Hasse diagram.
\[
\begin{tikzpicture}[xscale=2,yscale=.9,baseline=1.25cm]
\draw(0,1)node(1){$1$};
\draw(1,1)node(2){$2$};
\draw(0,2)node(3){$3$};
\draw(1,2)node(4){$4$};
\draw(0,3)node(5){$5$};
\draw(1,3)node(6){$6$};
\draw(0,4)node(7){$7$};
\draw(1,4)node(8){$8$};
\draw(6)--(4)--(1)--(3)--(5)--(7)--(4)--(2)--(5)--(8)--(6)--(3);
\end{tikzpicture}
\]

\begin{propn}\label{pnbc}
If $n\in\bbn$, then $\intn n$ is a \bc.
\end{propn}

\begin{pf}
$\intn n$ can be expressed as the union of two chains, namely the set of even integers in $\intn n$ and the set of odd integers in $\intn n$. This is the unique such expression, since if $\intn n$ is expressed in any other way as the union of two subsets, then one of these subsets contains two consecutive integers and so is not a chain.

Now suppose $\clsp$ is a proper refinement of $\cls$. Then there must be $i\in\{1,\dots,n-1\}$ such that either $i\clsp i+1$ or $i+1\clsp i$. Either way, we can find a new way to express $\intn n$ as the union of two $\clsp$-chains, namely
\[
\{\dots,i-4,i-2,i,i+1,i+3,i+5,\dots\}
\]
and
\[
\{\dots,i-3,i-1,i+2,i+4,\dots\}.
\]
This applies for any proper refinement $\clsp$, so $\intn n$ is a \bc.
\end{pf}

We now proceed to give several other characterisations of \bcs.

\section{Maximal elements}\label{maxlsec}

In this section we consider maximal elements in \bcs; this leads to some numerical results and a recursive construction for all \bcs. Several of the results in this section are due to Jeremy Rickard.

First we introduce some terminology. Say that a maximal element in a poset is \emph{supermaximal} if it lies above all the non-maximal elements.

\begin{lemmaciting}\textup{\textbf{(Rickard).}}\label{twomax}
Suppose $P$ is a \bc with $\card P\gs3$. Then $P$ contains exactly two maximal elements, exactly one of which is supermaximal.
\end{lemmaciting}

For example, looking at the four $5$-element \bcs depicted in \cref{bcsec}, we can see that each has exactly two maximal elements; in each diagram, the left-hand maximal element is supermaximal, but the right-hand one is not.

\begin{pf}
By definition $P$ is the union of two chains. No two maximal elements can lie in the same chain, so there can be at most two maximal elements in $P$. If there is only one maximal element $p$, then $p$ is the greatest element of $P$. But then there cannot be a unique way to express $P$ as the union of two chains, since if $B$ and $C$ are chains with $B\cup C=P$ and $p\in B$, then $B\setminus\{p\}$ and $C\cup\{p\}$ are also chains with union $P$.

So $P$ has exactly two maximal elements $p$ and $q$. We show by contradiction that exactly one of $p$ and $q$ is supermaximal.

If both $p$ and $q$ are supermaximal, then $P$ cannot be uniquely expressed as the union of two chains: indeed, if $P=B\cup C$ with $B,C$ chains, then without loss of generality we may assume that $p\in B$ and $q\in C$. But $p$ and $q$ lie above all other elements of $P$, so $B\setminus\{p\}\cup\{q\}$ and $C\setminus\{q\}\cup\{p\}$ are also chains whose union is $P$, so there is more than one way to write $P£$ as the union of two chains, a contradiction.

So suppose instead that neither $p$ nor $q$ is supermaximal. This means that there are $r,s\in P\setminus\{p,q\}$ such that $r\ncls q$ and $s\ncls p$. Now let $P=B\cup C$ be the unique decomposition of $P$ as the union of two chains. Then $p$ and $r$ lie in one of these chains, say $B$, and $q$ and $s$ lie in the other. Let $\clsp$ be the refinement of $\cls$ obtained by setting $r\clsp q$ and extending transitively; that is,
\[
a\clsp b\quad\Longleftrightarrow\quad a\cls b\ \text{ or }\ (a\cls r\text{ and }b=q).
\]
Since $\clsp$ is a proper refinement of $\cls$, the definition of a \bc means there are at least two ways to write $P$ as the union of two $\clsp$-chains. One of these is $P=B\cup C$; let $P=D\cup E$ be another. Note that $p$ and $q$ are $\clsp$-incomparable, and so are $p$ and $s$ and also $r$ and $s$. So $p$ and $r$ lie in one of the two chains, say $D$, while $q$ and $s$ lie in the other. But now we claim that $D$ and $E$ are also $\cls$-chains, which contradicts the assumption that $P$ can be uniquely written as the union of two $\cls$-chains. Given two elements of $D$, we can write them as $a\clsp b$, and since $b\neq q$, the definition of $\clsp$ means that $a\cls b$. So $D$ is a $\cls$-chain. Given two elements of $Q$, write them as $a\cls ^+ b$. Then the definition of $\clsp$ gives $a\cls b$ except possibly in the case where $a\cls r$ and $b=q$; so assume this is the case. Since $E$ is a $\clsp$-chain, $a$ must be $\clsp$-comparable with $s$, and hence $\cls$-comparable with $s$. We cannot have $s\cls a$, since then we would get $s\cls a\cls r\cls p$, contrary to assumptions. So instead $a\cls s$, and then $a\cls q=b$, as required. So in any case $a\cls b$, and $E$ is a $\cls$-chain.
\end{pf}

Now we give a reduction result which leads to a recursive construction of all \bcs.

\begin{propn}\textup{\textbf{(Rickard).}}\label{remsuper}
Suppose $P$ is a poset with $\card P\gs3$. Suppose $P$ has exactly two maximal elements $p$ and $q$, and that $p$ is supermaximal. Then $P$ is a \bc \iff $P\setminus\{p\}$ is a \bc.
\end{propn}

\begin{pf}
For this proof we write $P^-=P\setminus\{p\}$. First consider the case where $q$ is also supermaximal in $P$. Then neither $P$ nor $P^-$ is a \bc: $P^-$ fails to be a \bc because it has a greatest element $q$, so (as in the proof of \cref{twomax}) is not uniquely expressible as the union of two chains; $P$ fails to be a \bc by \cref{twomax} because it has two supermaximal elements.

So we can assume that $q$ is not supermaximal in $P$, and we consider decompositions of $P$ and $\pmi$ as unions of two chains. Suppose we can decompose $P$ as the union of two chains $B,C$, with $p\in B$ and $q\in C$. Then we can decompose $\pmi$ as the union of two chains $B\setminus\{p\}$ and $C$. Conversely, if we can decompose $\pmi$ as the union of two chains $D,E$ with $q\in E$, then $q\notin D$ (because $q$ is not supermaximal in $P$, and is therefore incomparable with some element of $P\setminus\{p\}$) so $P$ is the union of the two chains $D\cup\{p\}$ and $E$. Hence $P$ has a unique expression as the union of two chains \iff $\pmi$ does. So we can assume for the rest of the proof that $P$ and $\pmi$ each have a unique decomposition as the union of two chains.

Suppose $\pmi$ is not a \bc. Then there is a proper refinement $\clsp$ of $\cls$ on $\pmi$ such that $\pmi$ has a unique expression as the union of two $\clsp$-chains. If we extend $\clsp$ to all of $P$ by setting $r\clsp p$ for all $r\neq q$ and keeping $p,q$ incomparable, then $\clsp$ is a proper refinement of $\cls$ on $P$. Moreover, $P$ is uniquely expressible as the union of two $\clsp$-chains (since if there were more than one such expression, there would be more than one such expression for $\pmi$). So $P$ is not a \bc.

Conversely, suppose $P$ is not a \bc. Then there is a proper refinement $\clsp$ of $\cls$ such that $P$ is uniquely expressible as the union of two $\clsp$-chains. This means in particular that $P$ must have two $\clsp$-maximal elements (namely $p$ and $q$), so the restriction of $\clsp$ to $\pmi$ is a proper refinement of the restriction of $\cls$ to $\pmi$. Now $\pmi$ is uniquely expressible as the union of two $\clsp$-chains because $P$ is, so $\pmi$ is not a \bc.
\end{pf}

An immediate consequence is that we can determine all automorphisms of \bcs.

\begin{propn}\label{autotriv}
Suppose $P$ is a \bc with $\card P\gs3$. Then $P$ has no automorphisms except the identity.
\end{propn}

\begin{pf}
We use induction on $\card P$. The case where $\card P=3$ is trivial to check, so assume $\card P\gs4$. By \cref{twomax} $P$ has a unique supermaximal element $p$. Clearly $p$ must be fixed by any automorphism of $P$. Hence any automorphism of $P$ restricts to an automorphism of $P\setminus\{p\}$. By \cref{remsuper} $P\setminus\{p\}$ is a \bc, so by induction $P\setminus\{p\}$ has no non-trivial automorphisms. Hence neither does $P$.
\end{pf}

Rickard has pointed out that \cref{remsuper} gives a way to construct all \bcs recursively: starting from the unique $2$-element \bc, we repeatedly add new supermaximal elements; at each step, we simply have to choose which of the two existing maximal elements should remain maximal. For example, we may construct a $6$-element \bc by adding elements to form a sequence of \bcs as follows.
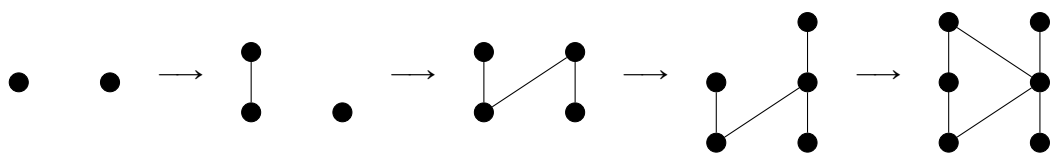
\[
\begin{tikzpicture}[scale=.8,baseline=0cm]
\nd{0,0}a
\nd{1.5,0}b
\end{tikzpicture}
\quad\longrightarrow\quad
\begin{tikzpicture}[scale=.8,baseline=.4cm]
\nd{0,0}a
\nd{1.5,0}b
\nd{0,1}c
\draw(c)--(a);
\end{tikzpicture}
\quad\longrightarrow\quad
\begin{tikzpicture}[scale=.8,baseline=.4cm]
\nd{0,0}a
\nd{1.5,0}b
\nd{0,1}c
\nd{1.5,1}d
\draw(c)--(a)--(d)--(b);
\end{tikzpicture}
\quad\longrightarrow\quad
\begin{tikzpicture}[scale=.8,baseline=.8cm]
\nd{0,0}a
\nd{1.5,0}b
\nd{0,1}c
\nd{1.5,1}d
\nd{1.5,2}f
\draw(f)--(d);
\draw(c)--(a)--(d)--(b);
\end{tikzpicture}
\quad\longrightarrow\quad
\begin{tikzpicture}[scale=.8,baseline=.8cm]
\nd{0,0}a
\nd{1.5,0}b
\nd{0,1}c
\nd{1.5,1}d
\nd{0,2}e
\nd{1.5,2}f
\draw(f)--(d)--(e)--(a)--(d)--(b);
\end{tikzpicture}
\]
We can make this precise in the following \lcnamecref{numberbcs}, in which we determine the number of $n$-element \bcs up to isomorphism.

\begin{propn}\label{numberbcs}
Suppose $n\gs3$. Then there are exactly $2^{n-3}$ \bcs of size $n$ up to isomorphism.
\end{propn}

\begin{pf}
We use induction on $n$, with the case $n=3$ being easy. Assuming $n\gs4$, \cref{remsuper} shows that we have a function $P\mapsto P^-$ from the set of $n$-element \bcs to the set of $(n-1)$-element \bcs, defined by removing the unique supermaximal element from a \bc. By the induction hypothesis it suffices to show that given an $(n-1)$-element \bc $Q$, there are exactly two \bcs $P$ up to isomorphism with $P^-=Q$. To reconstruct $P$ from $Q$, we just have to add a new element $p$ which is supermaximal in $P$, with $P$ having exactly one maximal element other than $p$. Since $Q$ has two maximal elements, we just choose which of these two elements will not lie below $p$ in $P$. So there are two possibilities for $P$; these are non-isomorphic, since any isomorphism from one to the other would restrict to a non-trivial automorphism of $Q$; but by \cref{autotriv} $Q$ has no non-trivial automorphisms.
\end{pf}

Next we show that the number of pairs of comparable elements in an $n$-element \bc depends only on $n$; this will be very useful in later sections.

\begin{propn}\label{noofpairs}
Suppose $P$ is a \bc with $\card P=n$. Then there are exactly $\binom n2+1$ pairs $(p,q)$ in $P$ with $p\cls q$, and hence exactly $n-1$ pairs $(p,q)$ with $p\ncls q\ncls p$.
\end{propn}

\begin{pf}
Suppose $n\gs3$. By \cref{twomax} $P$ has a supermaximal element $p$ and one other maximal element $q$, and by \cref{remsuper} $P\setminus\{p\}$ is a \bc. By induction there are exactly $\binom{n-1}2+1$ comparable pairs in $P\setminus\{p\}$; adding the supermaximal element $p$ adds $n-1$ comparable pairs (i.e.\ $r\cls p$ for every $r\neq q$), giving $\binom n2+1$ comparable pairs in $P$.
\end{pf}

We can also count ideals in \bcs.

\begin{propn}\label{bcideals}
Suppose $P$ is a \bc with $n$ elements, and $1\ls m<n$. Then $P$ has exactly two $m$-element ideals, and if $m\gs2$ then exactly one of these two ideals is a \bc.
\end{propn}

\begin{pf}
We use induction on $n$, with the case $n=2$ being trivial. Assuming $n\gs3$, let $p$ be the unique supermaximal element of $P$. Then by \cref{remsuper} $P^-=P\setminus\{p\}$ is a \bc. Since $p$ is supermaximal and there is only one other maximal element $q$, there are $n-2$ elements $r\in P$ with $r\cl p$. Hence $p$ cannot be contained in any ideal with fewer than $n-1$ elements. So if $m\ls n-2$, then any $m$-element ideal of $P$ is contained in $P^-$, and the result follows by the inductive hypothesis. So it remains to consider the case $m=n-1$. There are clearly two ideals with $n-1$ elements, namely $P^-$ and $P\setminus\{q\}$. Of these, $P^-$ is a \bc but $P\setminus\{q\}$ is not, because it contains only one maximal element, namely $p$.
\end{pf}

\begin{eg}
Consider the \bc $\intn n$ introduced in \cref{bcsec}. For $m\in\{1,\dots,n-1\}$, the two $m$-element ideals of $\intn n$ are $\{1,\dots,m\}$ and $\{1,\dots,m-1\}\cup\{m+1\}$. Of these, only the first is a \bc.
\end{eg}

In fact, there is a converse to \cref{bcideals}, which provides another characterisation of \bcs.

\begin{propn}\label{idealsbc}
Suppose $P$ is an $n$-element poset, and that $P$ has exactly two $m$-element ideals, for each $m\in\{1,\dots,n-1\}$. Then $P$ is a \bc.
\end{propn}

\begin{pf}
We use induction on $n$, and for the inductive step we assume $n\gs3$. Since $P$ has exactly two $(n-1)$-element ideals, it has exactly two maximal elements $p,q$. We claim that at least one of them must be supermaximal. If not, then there are elements $r,s\in P\setminus\{p,q\}$ with $r\ncls q$ and $s\ncls p$. If we take $r,s$ to be maximal with these properties, then $P$ has three $(n-2)$-element ideals
\[
P\setminus\{p,q\},\qquad P\setminus\{p,r\},\qquad P\setminus\{q,s\},
\]
a contradiction. So at least one of $p,q$, say $p$, is supermaximal. As observed in the proof of \cref{bcideals}, this means that every ideal with fewer than $n-1$ elements is contained in $P^-=P\setminus\{p\}$. So $P^-$ has exactly two $m$-element ideals, for each $m\in\{1,\dots,n-2\}$, so by induction $P^-$ is a \bc. Hence by \cref{remsuper} $P$ is a \bc too.
\end{pf}

\section{Two concrete constructions of \bcs}\label{consec}

Here we give a concrete construction of posets, which will turn out to be \bcs and in fact to give all \bcs up to isomorphism. This construction comes originally from the author's paper \cite{0hecke}.

First we set out some notation for binary sequences. A \emph{binary sequence} of length $n$ means a word $a=a_1\dots a_n$ in the alphabet $\{0,1\}$. Given such a sequence $a$, we define $\bar a = \bar a_1\dots\bar a_n$, where $\bar0=1$ and $\bar1=0$. We say that two binary sequences $b=b_1\dots b_r$ and $c=c_1\dots c_r$ are \emph{\ins} if for some $i,j$ we have $b_i>c_i$ and $b_j<c_j$, i.e.\ if the subsets of $\{1,\dots,r\}$ corresponding to $b$ and $c$ are incomparable in the containment order.

Now we can give our constructions. We fix a binary sequence $a=a_1\dots a_n$ of length $n\gs0$, and define $n+2$ sequences $a(0),\dots,a(n+1)$ of length $n+1$ as follows:
\[
\begin{array}{rc@{\,}c@{\,}c@{\,}c@{\,}c@{\,}c@{\,}c}
a(n+1)=&a_1&a_2&a_3&\dots&a_{n-1}&a_n&0\\
a(n)=&a_1&a_2&a_3&\dots&a_{n-1}&0&1\\
a(n-1)=&a_1&a_2&a_3&\dots&0&1&a_n\\
\hbox to 7pt{\hfil$\vdots$\hfil}\\
a(2)=&a_1&0&1&\dots&a_{n-2}&a_{n-1}&a_n\\
a(1)=&0&1&a_2&\dots&a_{n-2}&a_{n-1}&a_n\\
a(0)=&1&a_1&a_2&\dots&a_{n-2}&a_{n-1}&a_n.
\end{array}
\]

\begin{lemmac}{0hecke}{Proposition 6.1}\label{andist}
The sequences $a(0),\dots,a(n+1)$ are distinct.
\end{lemmac}

\begin{pf}
If the sequences $a(i)$ and $a(j)$ are equal for $i<j$, then
\[
1=a_{i+1}=a_{i+2}=\dots=a_{j-1}=0,
\]
a contradiction.
\end{pf}

In fact, there is another description of the set $\{a(0),\dots,a(n+1)\}$.

\begin{lemma}\label{insert}
$a(0),\dots,a(n+1)$ are precisely the sequences that can be obtained by inserting $0$ or $1$ at some point in $a$.
\end{lemma}

\begin{pf}
First we observe that each $a(i)$ is obtained by inserting a $0$ or a $1$ in $a$. For $i=0,n+1$ this is trivial, so suppose $1\ls i\ls n$. If $a_i=0$, then $a(i)$ is obtained from $a$ by inserting a $1$ immediately after $a_i$, while if $a_i=1$, then $a(i)$ is obtained from $a$ by inserting a $0$ immediately before $a_i$.

Conversely, suppose $b$ is the sequence obtained by inserting $x\in\{0,1\}$ immediately after position $i$. We assume $x=0$, as the case $x=1$ is similar. Let $j>i$ be maximal such that $a_{i+1}=\dots=a_{j-1}=0$. Then $b=a(j)$; indeed, both sequences equal
\[
a_1\dots a_i00\dots01a_{j+1}\dots a_n
\]
where there are $j-i$ $0$s, and the $1$ should be omitted if $j=n+1$.
\end{pf}

The description of the set $\{a(0),\dots,a(n+1)\}$ given by \cref{insert} is simpler than the original definition. However, the ordering $a(0),\dots,a(n+1)$ is significant for defining our partial order, which we do next.

Set $\seqp a=\{a(0),\dots,a(n+1)\}$ and define a binary relation $\cls$ on $\seqp a$ by setting $a(i)\cl a(j)$ whenever $i<j$ and $a(i)$ and $a(j)$ are \ins.

\begin{eg}
Take $a=1101$. Then we have
\begin{align*}
a(5)&=11010\\
a(4)&=11001\\
a(3)&=11011\\
a(2)&=10101\\
a(1)&=01101\\
a(0)&=11101
\end{align*}
and the relation $\cls$ is a partial order, with the following Hasse diagram.
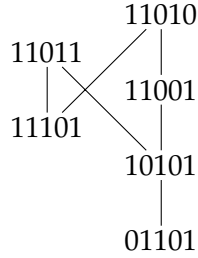
\[
\begin{tikzpicture}[scale=1,xscale=-1,every node/.style={inner sep=1.5pt}]
\draw(0,0)node(a){$01101$};
\draw(0,1)node(b){$10101$};
\draw(0,2)node(c){$11001$};
\draw(0,3)node(d){$11010$};
\draw(1.5,1.5)node(e){$11101$};
\draw(1.5,2.5)node(f){$11011$};
\draw(a)--(b)--(c)--(d)--(e)--(f)--(b);
\end{tikzpicture}
\]
\end{eg}

It is not obvious from the definition that $\cls$ is a partial order in general -- in particular, that it is transitive. In fact this was shown in \cite{0hecke}, where $\seqp a$ is shown to be the ``structure poset'' of a multiplicity-free module for a $0$-Hecke algebra. We will show it in a different (and purely combinatorial) way by giving a different definition of $\cls$.

Keeping the binary sequence $a$ fixed, define a second list of binary sequences $a[0],\dots,a[n+1]$ by
\[
\begin{array}{rc@{\,}c@{\,}c@{\,}c@{\,}c@{\,}c@{\,}c}
a[n+1]=&a_1&a_2&a_3&\dots&a_{n-1}&a_n&1\\
a[n]=&a_1&a_2&a_3&\dots&a_{n-1}&1&0\\
a[n-1]=&a_1&a_2&a_3&\dots&1&0&a_n\\
\vdots\\
a[2]=&a_1&1&0&\dots&a_{n-2}&a_{n-1}&a_n\\
a[1]=&1&0&a_2&\dots&a_{n-2}&a_{n-1}&a_n\\
a[0]=&0&a_1&a_2&\dots&a_{n-2}&a_{n-1}&a_n.
\end{array}
\]

\begin{lemma}\label{samelist}
$\seqp a=\{a[0],\dots,a[n+1]\}$.
\end{lemma}

\begin{pf}
\cref{insert} says that $\seqp a=\{a(0),\dots,a(n+1)\}$ is the set of all sequences that can be obtained by inserting a symbol in $a$. But the latter description is symmetric in $0$ and $1$, and therefore the same applies to the set $\{a[0],\dots,a[n+1]\}$.
\end{pf}

As a consequence of \cref{samelist}, we have a permutation $w_a$ of $\{0,\dots,n+1\}$ defined by $a(i)=a[w_a(i)]$. This allows us to give a second description of the relation $\cls$. First we give a simple description of the permutation $w_a$; we leave the proof as an exercise.

\begin{lemma}\label{simplewa}
Suppose $i\in\{0,\dots,n+1\}$.

\begin{itemize}
\item
If $i=0$ or $a_i=0$, then $w_a(i)$ is the smallest $k>i$ such that $k=n+1$ or $a_k=0$.
\item
If $i=n+1$ or $a_i=1$, then $w_a(i)$ is the largest $k<i$ such that $k=0$ or $a_k=1$.
\end{itemize}
\end{lemma}

\begin{propn}\label{posequal}
Given $i,j\in\{0,\dots,n+1\}$ we have $a(i)\cls a(j)$ \iff $i\ls j$ and $w_a(i)\ls w_a(j)$. Hence $\cls$ is a partial order on $\seqp a$.
\end{propn}

\begin{pf}
We need to show that for $0\ls i<j\ls n+1$ the sequences $a(i)$ and $a(j)$ are \ins \iff $w_a(i)<w_a(j)$. We consider several cases; in the remainder of this proof we read $a_0$ as $0$ and $a_{n+1}$ as $1$.
\begin{enumerate}[ref=(\arabic{enumi})]
\item\label{case1}
Suppose $a_i=a_j$. Then $a(i)$ and $a(j)$ are certainly \ins because they are different but have the same sum. Moreover, $w_a(i)<w_a(j)$, by \cref{simplewa}.
\item
Suppose $a_i=1$ and $a_j=0$. Then $a(i)$ and $a(j)$ are \ins, since
\[
a(i)_i=a(j)_k=0,\qquad a(i)_k=a(j)_i=1,
\]
where $k>i$ is minimal such that $a_k=0$. Furthermore, $w_a(i)<i<j<w_a(j)$ by \cref{simplewa}.
\item
Finally suppose $a_i=0$ and $a_j=1$.
\begin{itemize}
\item
If $w_a(j)<i$, then certainly $w_a(i)>w_a(j)$. And in this case $a_{i+1}=\dots=a_{j-1}=0$ by \cref{simplewa}, so that $a(i)$ and $a(j)$ differ only in the $(i+1)$th position and are not \ins.
\item
Similarly if $w_a(i)>j$, then $w_a(i)>w_a(j)$ and $a(i),a(j)$ are not \ins.
\item
If $w_a(i)=w_a(j)+1$, then the sequence $a_i a_{i+1}\dots a_j$ has the form $011\dots1100\dots001$. So $a(i)$ and $a(j)$ are not \ins: they differ only in the $w_a(i)$th position.
\item
The remaining possibility is that $i<w_a(i)<w_a(j)<j$. In this case let $k=w_a(i)$, and let $l>k$ be minimal such that $a_l=1$. Then by assumption $l<j$, so
\[
a(i)_k=a(j)_l=1,\qquad a(i)_l=a(j)_k=0,
\]
and $a(i)$ and $a(j)$ are \ins.\qedhere
\end{itemize}
\end{enumerate}
\end{pf}

\cref{posequal} shows that $\cls$ is the intersection of the two total orders
\[
a(0)<\dots<a(n+1)\qquad\text{and}\qquad
a[0]<\dots<a[n+1],
\]
so the poset $\seqp a$ has dimension $2$ (in general, the \emph{dimension} of a poset $P$ is the smallest $d$ such that the partial order on $P$ can expressed as the intersection of $d$ total orders). In fact, this follows from a more general result of Pretzel \cite[Theorem 1]{prt}, which says that the dimension of a poset is no more than its width, i.e.\ the size of its largest antichain; clearly in a \bc the largest antichain has size $2$.

In \cref{coxsec} we will examine the permutations $w_a$ in more detail.

\begin{cory}\label{01symm}
There is an isomorphism
\begin{align*}
\phi:\seqp a&\longrightarrow \seqp{\bar a}\\
a(i)&\longmapsto\widebar{a(i)}.
\end{align*}
\end{cory}

\begin{pf}
Notice that $\widebar{a(i)}=\bar a[i]$ and $\widebar{a[i]}=\bar a(i)$ for each $i$, so $\phi$ is certainly a bijection from $\seqp a$ to $\seqp{\bar a}$, and \cref{posequal} guarantees that $\phi(a(i))\cls\phi(a(j))$ \iff $a(i)\cls a(j)$.
\end{pf}

We note in passing that this \lcnamecref{01symm} yields a third way to describe the partial order on $\seqp a$: we have $a[i]\cl a[j]$ \iff $i<j$ and $a[i]$ and $a[j]$ are \ins.

Now we prove the main result of this section.

\begin{thm}\label{sameposets}
Suppose $n\gs0$.
\begin{enumerate}[ref=(\arabic*)]
\item\label{seqpbc}
If $a$ is a binary sequence of length $n$, then $\seqp a$ is a \bc.
\item
If $P$ is an $(n+2)$-element \bc, then $P\cong\seqp a$ for some binary sequence $a$.
\item
If $a,b$ are binary sequences with $\seqp a\cong\seqp b$, then $b=a$ or $b=\bar a$.
\end{enumerate}
\end{thm}

\begin{pfenum}
\item
We use induction on $n$, and \cref{remsuper}. The case $n=0$ is trivial, so take $n\gs1$. Using \cref{01symm} we can assume that $a_n=1$. First we claim that $\seqp a$ has exactly two maximal elements $a(n+1)$ and $a[n+1]$, with $a(n+1)$ being supermaximal. \cref{posequal} certainly shows that $a(n+1)$ and $a[n+1]$ are maximal. Furthermore, the assumption $a_n=1$ means that $a(n+1)=a[n]$. So (again using \cref{posequal}) $a(n+1)$ lies above every element of $\seqp a$ except $a[n+1]$. So $a(n+1)$ is supermaximal, and there are no maximal elements other than $a(n+1)$ and $a[n+1]$.

Now define $a^-$ to be the binary sequence $a_1\dots a_{n-1}$. Let $\seqp a^-=\{a(0),\dots,a(n)\}$, with the partial order $\cls$ induced from $\seqp a$. Then we observe that $\seqp a^-$ is isomorphic to $\seqp{a^-}$: each of the sequences $a(0),\dots,a(n)$ ends in a $1$, and deleting this $1$ from each sequence yields the sequences $a^-(0),\dots,a^-(n)$, preserving \insness. So by induction $\seqp a^-$ is a \bc, and hence by \cref{remsuper} so is $\seqp a$.
\item
Again we use induction on $n$, and the case $n=0$ is trivial. By \cref{twomax} $P$ has one supermaximal element $p$, and one other maximal element $q$. Let $\pmi=P\setminus\{p\}$. Then by \cref{remsuper} $\pmi$ is a \bc, so by induction there is an isomorphism $\pmi\cong\seqp{a^-}$ for some binary sequence $a^-$ of length $n-1$. Under this isomorphism, the two maximal elements of $\pmi$ (one of which is $q$) map to $a^-(n)$ and $a^-[n]$. We suppose $q$ maps to $a^-[n]$ (the other case is similar). Define a sequence $a$ by adding a $1$ at the end of $a^-$; then we claim that $P\cong \seqp a$. As explained in the proof of \ref{seqpbc}, $a(n+1)$ is a supermaximal element of $\seqp a$, and $\seqp a\setminus\{a(n+1)\}\cong \seqp{a^-}$, and so there is an isomorphism $\theta:\seqp a\setminus\{a(n+1)\}\to\pmi$, with $\theta(a[n+1])=q$. Since $a(n+1)$ is incomparable with $a[n+1]$ in $\seqp a$ and lies above every other element, and $p$ is incomparable with $q$ in $P$ but above every other element, we can extend $\theta$ by setting $\theta(a(n+1))=p$ to give an isomorphism from $\seqp a$ to $P$.
\item
Again, we use induction on $n$, and assume here that $n\gs2$. Let $\phi:\seqp a\to\seqp b$ be an isomorphism; in fact, by \cref{autotriv} $\phi$ is the unique isomorphism from $\seqp a$ to $\seqp b$. Then $\phi$ must map the supermaximal element of $\seqp a$ to the supermaximal element of $\seqp b$, and the other maximal element of $\seqp a$ to the other maximal element of $\seqp b$. Arguing as in the proof of \ref{seqpbc} above, we see that the supermaximal element of $\seqp a$ is $a_1\dots a_n\widebar{a_n}$, and the other maximal element is $a_1\dots a_na_n$. A similar statement holds for $\seqp b$, so $\phi$ maps $a_1\dots a_na_n$ to $b_1\dots b_nb_n$.

Let $\seqp a^-$ be the \bc obtained by deleting the supermaximal element of $\seqp a$, and define $\seqp b^-$ similarly. Then (as in the proof of \ref{seqpbc}, interchanging $0$ and $1$ if necessary) $\seqp a^-\cong\seqp{a^-}$, with the unique isomorphism $\theta:\seqp a^-\to\seqp{a^-}$ defined by deleting the last digit $a_n$ from the end of each element of $\seqp a^-$. Similarly, there is a unique isomorphism $\kappa:\seqp b^-\to\seqp{b^-}$, given by deleting the digit $b_n$ at the end of each element of $\seqp b^-$.

Now we have an isomorphism $\kappa\circ\phi\circ\theta^{-1}:\seqp{a^-}\to\seqp{b^-}$ which maps $a$ to $b$. By induction we see that $b^-$ equals either $a^-$ or $\widebar{a^-}$. But the unique isomorphism from $\seqp{a^-}$ to $\seqp{a^-}$ maps $a$ to $a$, while the unique isomorphism from $\seqp{a^-}$ to $\seqp{\widebar{a^-}}$ maps $a$ to $\bar a$ by \cref{01symm}. So $b$ must equal either $a$ or $\bar a$.
\end{pfenum}

\cref{sameposets} allows us to label \bcs in a canonical way: each \bc can be written as $\seqp a$ for some binary sequence $a$ which is unique up to replacing $a$ with $\bar a$. Note that given \cref{numberbcs,01symm}, any two parts of \cref{sameposets} imply the other part. But proving all three parts directly yields a new proof of \cref{numberbcs}.

\smallskip
Now we consider duality. Let us write $P\dual$ for the poset dual to $P$. It is obvious from the definition in \cref{bcsec} that the dual of a \bc is again a \bc. Using the labelling for \bcs provided by \cref{sameposets}, we can be more specific.

\begin{propn}\label{dual}
Suppose $a$ is a binary sequence of length $n$, and let $\rev a$ be the binary sequence $a_n\dots a_1$. Then $\seqp a\dual\cong\seqp{\rev a}$.
\end{propn}

\begin{pf}
This follows from \cref{posequal}, which says that $\cls$ is the intersection of the total orders
\[
a(0)<\dots<a(n+1)\qquad\text{and}\qquad
a[0]<\dots<a[n+1]
\]
on $\seqp a$. Hence the partial order on $\seqp a\dual$ is the intersection of the total orders
\[
a(n+1)<\dots<a(0)\qquad\text{and}\qquad
a[n+1]<\dots<a[0].
\]
But now observe that $(\rev a)(i)=\rev{(a[n+1-i])}$ and $(\rev a)[i]=\rev{(a(n+1-i))}$ for each $i$. So we have an isomorphism from $\seqp a\dual$ to $\seqp{\rev a}$ given by mapping $a(i)\mapsto\rev{(a(i))}$ for each~$i$.
\end{pf}

As a consequence of \cref{dual}, we see that a \bc $P=\seqp a$ is self-dual \iff $\rev a$ equals either $a$ or $\bar a$. Clearly a self-duality preserves the unique decomposition of $P$ as a union of two chains; one can show that the self-duality preserves the two chains if $\rev a=a$, and interchanges them if $\rev a=\bar a$.

\section{Splicing \bcs}\label{splicesec}

In this section we give a simple way of joining two \bcs together to create a larger one. This yields a further construction of all \bcs, starting from a family of ``indecomposable'' \bcs.

We start with a more general definition. Suppose $P$ and $Q$ are posets, and that $P$ has exactly two maximal elements $p_0,p_1$, with only $p_0$ being supermaximal, and suppose $Q$ has exactly two minimal elements $q_0,q_1$ with only $q_0$ being superminimal (i.e.\ lying below every non-minimal element). We define a new poset $P\splice Q$ called the \emph{splice} of $P$ and $Q$. Informally, this is defined by placing $Q$ above $P$, and identifying $p_0$ with $q_0$ and $p_1$ with $q_1$. Formally, we start by defining a partial order on the disjoint union $P\sqcup Q$ via the following rules:
\begin{itemize}
\item
if $a,b\in P$, then $a\cls b$ in $P\sqcup Q$ \iff $a\cls b$ in $P$;
\item
if $a,b\in Q$, then $a\cls b$ in $P\sqcup Q$ \iff $a\cls b$ in $Q$;
\item
if $a\in P$ and $b\in Q$, then $a\cls b$ \iff for some $i\in\{0,1\}$ we have $a\cls p_i$ in $P$ and $q_i\cls b$ in~$Q$;
\item
if $a\in Q$ and $b\in P$, then $a\ncls b$.
\end{itemize}
Now define $P\splice Q$ to be the quotient poset obtained by identifying $p_0$ with $q_0$ and $p_1$ with $q_1$. Note that $P$ and $Q$ are then naturally subposets of $P\splice Q$.

Now we consider the special case of \bcs. Recall that a \bc of size at least $3$ has exactly two maximal elements, with exactly one being supermaximal. Dually, a \bc of size at least $3$ has exactly two minimal elements, with exactly one being superminimal. Hence the splice of two \bcs each with at least three elements is well-defined. In fact we can easily extend this definition to include the case where either or both of the posets is the unique $2$-element \bc $\seqp\varnothing$, even though both elements of $\seqp\varnothing$ are both supermaximal and superminimal. We then get $P\splice\seqp\varnothing=P$ and $\seqp\varnothing\splice Q=Q$.

\begin{propn}\label{splicebc}
If $P$ and $Q$ are \bcs, then so is $P\splice Q$.
\end{propn}

\begin{pf}
Since $P$ and $Q$ can each be expressed as the union of two chains, so can $P\splice Q$: we take the union of the chain in $P$ containing $p_0$ and the chain in $Q$ containing $q_0$ to give one chain, and do the same for $p_1,q_1$ to give the other. Moreover, this is the unique way to express $P\splice Q$ as the union of two chains, since any other such expression would restrict to give a new expression for either $P$ or $Q$ as the union of two chains.

Now suppose $\clsp$ is a proper refinement of $\cls$ on $P\splice Q$. Note that for every $p\in(P\splice Q)\setminus Q$ and $q\in(P\splice Q)\setminus P$ we already have $p\cls q$, since $p\cls p_0$ in $P$ and $q_0\cls q$ in $Q$. So in order for $\clsp$ to be a proper refinement of $\cls$ on $P\splice Q$, we must have $a\clsp b$ but $a\ncls b$ either for some $a,b\in P$ or for some $a,b\in Q$. We assume the former case; then the restriction of $\clsp$ to $P$ is a proper refinement of $\cls$ on $P$. Since $P$ is a \bc, this means that there are at least two different ways to express $P$ as the union of two $\clsp$-chains. But now (via the construction in the first paragraph of the proof) there are at least two ways to express $P\splice Q$ as the union of two $\clsp$-chains.
\end{pf}

\begin{eg}
Take $P=Q=\seqp{001}$. Then the Hasse diagram of $P$ with the elements $p_0,p_1,q_0,q_1$ marked is as follows.
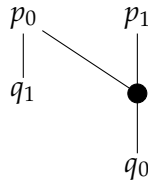
\[
\begin{tikzpicture}[scale=1,yscale=-1,every node/.style={inner sep=1.5pt}]
\draw(0,0)node(a){$p_0$};
\draw(0,1)node(b){$q_1$};
\draw(1.5,0)node(c){$p_1$};
\nd{1.5,1}d
\draw(1.5,2)node(e){$q_0$};
\draw(b)--(a)--(d);
\draw(c)--(e);
\end{tikzpicture}
\]
Hence the splice $P\splice Q$ is given by the following diagram.
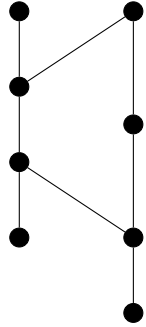
\[
\begin{tikzpicture}[yscale=1,xscale=1.5]
\nd{0,0}a
\nd{0,1}b
\nd{0,2}c
\nd{0,3}d
\nd{1,-1}e
\nd{1,0}f
\nd{1,1.5}g
\nd{1,3}h
\draw(a)--(b)--(c)--(d);
\draw(e)--(f)--(g)--(h)--(c);
\draw(f)--(b);
\end{tikzpicture}
\]
\end{eg}

We can make \cref{splicebc} more explicit using the labelling for \bcs introduced in \cref{consec}. Given binary sequences $a,b$, write $a|b$ for their concatenation.

\begin{propn}\label{finitesplice}
Suppose $a=a_1\dots a_r$ and $b=b_1\dots b_s$ are binary sequences, and define
\[
c=
\begin{cases}
a|b&\text{if $a_r=b_1$}\\
a|\bar b&\text{if $a_r\neq b_1$}.
\end{cases}
\]
Then $\seqp a\splice \seqp b\cong \seqp c$.
\end{propn}

(Note that we include the case where $r$ or $s$ equals $0$; in this case the condition $a_r=b_1$ doesn't make sense, but the conclusion is trivially true regardless.)

\begin{pf}
By replacing $b$ with $\bar b$ if necessary and using the fact that $\seqp b\cong\seqp{\bar b}$, we can assume that $a_r=b_1$. In fact, by \cref{01symm} we can assume that $a_r=0=b_1$. Then the supermaximal element of $\seqp a$ is $a|1=a[r+1]$, with the other maximal element being $a|0=a(r+1)$. Similarly, the superminimal element of $\seqp b$ is $1|b=b(0)$, and the other minimal element is $0|b=b[0]$.

We define an injective function
\begin{align*}
\phi_b:\seqp a&\longrightarrow\seqp{a|b}\\
a[i]&\longmapsto (a|b)[i].
\end{align*}
This function can be more simply described as mapping $p\mapsto p|b$ for each $p\in \seqp a$, which shows that it is order-preserving.
\endgraf
We also define an injective function
\begin{align*}
\phi^a:\seqp b&\longrightarrow\seqp{a|b}\\
b(i)&\longmapsto (a|b)(r+i).
\end{align*}
This function can be more simply described as mapping $q\longmapsto a|q$, which shows that it too is order-preserving.

The images of $\phi_b$ and $\phi^a$ intersect in the two points
\begin{alignat*}2
\phi_b(a|0)&=\,&a|0|b&=\phi^a(0|b)\\
\intertext{and}
\phi_b(a|1)&=\,&a|1|b&=\phi^a(1|b).
\end{alignat*}

Since $a|0$ is supermaximal in $\seqp a$ and $0|b$ is superminimal in $\seqp b$, we get
\[
p\cls a|0|b\cls q
\]
for all $p\in\im\phi_b\setminus\{a|1|b\}$ and $q\in\im\phi^a\setminus\{a|1|b\}$. Putting this together with the partial orders on $\im\phi_b$ and $\im\phi^a$, we see that $\seqp{a|b}\cong \seqp a\splice \seqp b$.
\end{pf}

As a consequence of \cref{finitesplice}, we can see that whenever a binary sequence $a$ contains two consecutive equal entries, $\seqp a$ can be expressed as a splice of two smaller \bcs. This means that we can construct all \bcs from the \bcs of the form $\seqp{1010\dots}$. In fact, we have already seen these \bcs constructed in a very simple way: they are the posets $\intn n$ defined in \cref{bcsec}, as we now show.

\begin{propn}\label{altbc}
Suppose $a$ is the alternating binary sequence $1010\dots$ of length $n$. Then $\seqp a\cong\intn{n+2}$.
\end{propn}

\begin{pf}
Suppose $0\ls i<j\ls n+1$. We can check that $a(i)$ and $a(j)$ are \ins unless (and only unless) one of the following happens:
\begin{itemize}
\item
$i$ is even and $j=i+1$;
\item
$i$ is even and $j=i+3$;
\item
$n$ is odd, $i=n-1$ and $j=n+1$.
\end{itemize}
Hence we can define an isomorphism from $\seqp a$ to $\intn{n+2}$ by mapping
\[
a(i)\longmapsto
\begin{cases}
i&\text{if $i$ is odd}\\
i+2&\text{if $i$ is even and $i\neq n+1$}\\
i+1&\text{if $i$ is even and $i=n+1$}.
\end{cases}\qedhere
\]
\end{pf}

As a consequence, we deduce the following.

\begin{propn}\label{splicepn}
Suppose $P$ is a \bc with $n$ elements, where $n\gs3$. Then $P$ is isomorphic to $\intn{n_1}\splice\cdots\splice \intn{n_s}$ for a unique choice of $n_1,\dots,n_s\gs3$ with $n_1+\dots+n_s=n+2s-2$.
\end{propn}

\begin{pf}
The existence of such a decomposition comes from \cref{finitesplice,altbc}: writing $P\cong \seqp a$ for a binary sequence $a$, we break $a$ into subsequences $a^{(1)},\dots,a^{(s)}$ where each $a^{(i)}$ is an alternating sequence $\dots0101\dots$, and the last term of $a^{(i)}$ equals the first term of $a^{(i+1)}$ for each $i$. Then $n_1,\dots,n_s$ are just the lengths of the sequences $a^{(1)},\dots,a^{(s)}$.

For the uniqueness, it suffices to use \cref{numberbcs} and to count the number of possible expressions $n_1+\dots+n_s=n+2s-2$, which is a simple exercise.
\end{pf}

\begin{eg}
The \bc $\seqp{011011}$ is isomorphic to $\intn4\splice\intn5\splice\intn3$, as we see from the following Hasse diagrams.
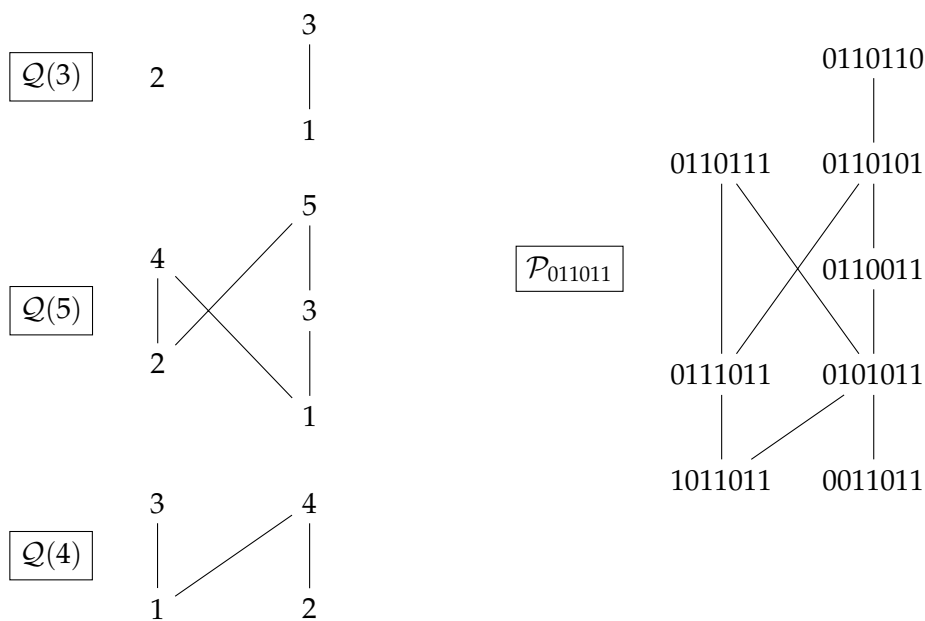
\[
\begin{array}{c}
\begin{tikzpicture}[xscale=2,yscale=.7,baseline=1.25cm]
\draw(-.7,2)node[draw]{$\intn3$};
\draw(1,1)node(1){$1$};
\draw(0,2)node(2){$2$};
\draw(1,3)node(3){$3$};
\draw(3)--(1);
\end{tikzpicture}
\\\\
\begin{tikzpicture}[xscale=-2,yscale=.7]
\draw(1.7,3)node[draw]{$\intn5$};
\draw(0,1)node(1){$1$};
\draw(1,2)node(2){$2$};
\draw(0,3)node(3){$3$};
\draw(1,4)node(4){$4$};
\draw(0,5)node(5){$5$};
\draw(3)--(1)--(4)--(2)--(5)--(3);
\end{tikzpicture}
\\\\
\begin{tikzpicture}[xscale=2,yscale=.7]
\draw(-.7,2)node[draw]{$\intn4$};
\draw(0,1)node(1){$1$};
\draw(1,1)node(2){$2$};
\draw(0,3)node(3){$3$};
\draw(1,3)node(4){$4$};
\draw(3)--(1)--(4)--(2);
\end{tikzpicture}
\end{array}
\qquad\qquad\qquad
\begin{tikzpicture}[scale=1,xscale=-2,yscale=1.4,baseline=2cm,every node/.style={fill=white}]
\draw(2,2)node[draw]{$\seqp{011011}$};
\draw(0,0)node(f){$0011011$};
\draw(0,1)node(e){$0101011$};
\draw(0,2)node(c){$0110011$};
\draw(0,3)node(b){$0110101$};
\draw(0,4)node(a){$0110110$};
\draw(1,0)node(h){$1011011$};
\draw(1,1)node(g){$0111011$};
\draw(1,3)node(d){$0110111$};
\draw(f)--(e)--(c)--(b)--(a);
\draw(b)--(g)--(d)--(e)--(h)--(g);
\end{tikzpicture}
\]
\end{eg}

\section{Graphs}\label{graphsec}

In this section we consider an analogue of \bcs in the context of graphs. We consider only finite graphs. Let's say that a graph $G=(V,E)$ is a \emph{\bicl} if $V$ can be uniquely expressed as the union of two cliques (i.e.\ complete subgraphs) and $G$ is maximal with this property: adding any edge to $E$ breaks the uniqueness.

In fact, \bicls are easy to understand, by looking at complements. A graph is a \bicl \iff its complement is bipartite with a unique bipartition, and is minimal with this property (i.e.\ removing any edge breaks the uniqueness). Being bipartite with a unique bipartiition is the same as being connected and bipartite, and a minimal connected (bipartite) graph is a tree, so a \bicl is simply the complement of a tree.

However, \bicls do have a direct relationship with \bcs. The \emph{comparability graph} of a poset $P$ has vertex set $P$, with an edge from $p$ to $q$ \iff either $p\cl q$ or $q\cl p$. Our main result here is the following, which gives yet another characterisation of \bcs. The author is grateful to an anonymous referee for suggesting the proof of this result, which is simpler than the proof in the original version of this paper.

\begin{propn}\label{bichainbiclique}
Suppose $P$ is a poset. Then $P$ is a \bc \iff its comparability graph is a \bicl.
\end{propn}

\begin{pf}
Let $G$ be the \emph{in}comparability graph of $P$, in which there is an edge from $p$ to $q$ \iff $p\ncls q\ncls p$. From the above discussion, we need to show that $P$ is a \bc \iff $G$ is a tree. Let $n=\card P$, and proceed by induction on $n$, with the case $n=2$ being trivial.

First consider the number of maximal elements of $P$. If $P$ has only one maximal element $p$, then by \cref{twomax} $P$ is not a \bc and $G$ is not a tree, because $p$ is an isolated vertex of $G$. If $P$ has three (or more) maximal elements, then by \cref{twomax} $P$ is not a \bc, and $G$ is not a tree, since it contains a triangle.

So we can assume that $P$ has exactly two maximal elements. Furthermore, we can assume that one of these maximal elements is supermaximal: if not, then by \cref{twomax} $P$ is not a \bc, and $G$ contains a $4$-cycle, so is not a tree. So let $p$ be a supermaximal element of $P$, let $P^-=P\setminus\{p\}$, and let $G^-$ be the incomparability graph of $P^-$. Since $p$ is supermaximal, there is only one element of $P$ incomparable with $p$ in $P$, so $G$ is obtained from $G^-$ by adding a leaf. So $G$ is a tree \iff $G^-$ is a tree, and by \cref{remsuper} $P$ is a \bc \iff $P^-$ is a \bc, so the result follows by induction.
\end{pf}

One can ask which trees can occur as the incomparability graphs of posets (and hence of \bcs). The general question of which graphs are incomparability graphs of posets was answered completely by Gallai \cite{gallai}. His theorem gives an explicit list of graphs, and says that a graph is the incomparability graph of a poset \iff it has no induced subgraph in the given list. We are interested in the special case of trees; for this, we just check that the only forest in Gallai's list is the \emph{triskelion}
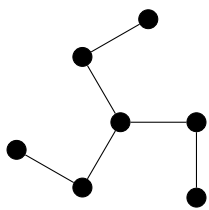
\[
\begin{tikzpicture}[scale=1]
\nd{0,0}a
\nd{1,0}b
\nd{1,-1}c
\nd{120:1}d
\nd{120:1)++(30:1}e
\nd{240:1}f
\nd{240:1)++(150:1}g
\draw(c)--(b)--(a)--(d)--(e);
\draw(a)--(f)--(g);
\end{tikzpicture}
\]
which leads to the following result.

\begin{thmciting}{\textbf{\textup{(Gallai).}}}\label{gallaitree}
Suppose $G$ is a tree. Then $G$ is the incomparability graph of a poset \iff $G$ is a \emph{caterpillar}, i.e.\ if the non-leaves in $G$ form a path.
\end{thmciting}

\cref{bichainbiclique,gallaitree} show that the incomparability graph of a \bc is a caterpillar, and that every caterpillar arises as the incomparability graph of a \bc. We now show that this \bc is unique up to duality.

\begin{propn}\label{catoccurs}
Suppose $G$ is a caterpillar. Then $G$ is the incomparability graph of a \bc $P$, which is unique up to duality.
\end{propn}

\begin{pf}
We use a counting argument. Sending a poset to its incomparability graph defines a function from (isomorphism classes of) \bcs to caterpillars, and \cref{bichainbiclique,gallaitree} show that this function is surjective. Moreover, dual \bcs obviously map to the same caterpillar. So it suffices to show that the number of dual pairs of $n$-element \bcs equals the number of isomorphism types of caterpillars with $n$ vertices. \cref{dual} and the discussion preceding \cref{numberbcs} show that the former number is the number of binary sequences of length $n-2$ modulo reversal and modulo interchanging $0$s and $1$s; an easy exercise shows that this number is $2^{n-4}+2^{\lfloor(n-4)/2\rfloor}$. \cite[Theorem 2.1]{hasc} shows that this is also the number of isomorphism types of caterpillars on $n$ vertices.
\end{pf}

By way of example, we consider two special cases of caterpillars. The first example is the \emph{star}, where one vertex is attached to all the others. This is the incomparability graph of the \bc consisting of a chain and an isolated point. The second example is the \emph{path} with $n$ vertices: this is the incomparability graph of the \bc $\intn n$ introduced in \cref{bcsec}.

\section{Coxeter elements in the symmetric group}\label{coxsec}

In this section we give yet another construction of \bcs, relating to permutations. Recall from \cref{consec} the second construction of the poset $\seqp a$: we define a permutation $w_a$ of $\{0,\dots,n+1\}$ by $a(i)=a[w_a(i)]$; then we have $a(i)\cls a(j)$ in $\seqp a$ if $i\ls j$ and $w_a(i)\ls w_a(j)$. Our focus in this section is on the permutations $w_a$: we consider the analogous definition for an arbitrary permutation, and determine exactly which permutations yield \bcs.

We recall some basic Coxeter theory of the symmetric group (for more details see the book by Humphreys \cite{hum}). Let $\sss n$ denote the group of all permutations of $\{1,\dots,n\}$. We write $s_1,\dots,s_{n-1}$ for the \emph{Coxeter generators} of $\sss n$; here $s_i$ is the transposition $(i,\ipo)$. Any permutation $w\in\sss n$ can be written in the form $s_{i_1}\dots s_{i_l}$ for some $i_1,\dots,i_l$; the smallest $l$ for which we can do this is called the \emph{length} $l(w)$. The length of $w$ can be also be defined as the number of \emph{inversions} of $w$, i.e.\ pairs $i<j$ such that $w(i)>w(j)$.

A \emph{Coxeter element} of $\sss n$ is an element which can be written as a product of the Coxeter generators, each appearing exactly once, in some order. Alternatively, a Coxeter element may be described as an element $w\in\sss n$ such that for every $r\neq1,n$, one of $w(r)$, $w^{-1}(r)$ is greater than $r$ and the other is less than $r$.

\begin{eg}
The Coxeter elements in $\sss4$ are the $4$-cycles
\[
(1,2,3,4),\qquad(1,3,4,2),\qquad(1,2,4,3),\qquad(1,4,3,2).
\]
\end{eg}

Coxeter elements are defined for all finite Coxeter groups, and play an important role in Coxeter theory. For the symmetric group, we will need the following result.

\begin{propn}\label{coxprops}
Suppose $n\gs2$. Then there are exactly $2^{n-2}$ Coxeter elements in $\sss n$, and they all have length $n-1$ and order $n$. In particular, if $n\gs3$, then no Coxeter element of $\sss n$ is self-inverse.
\end{propn}

The number of Coxeter elements in $\sss n$ is a special case of a theorem of Shi \cite[Theorem 1.5]{shi}. The order of the Coxeter element $s_1\dots s_{n-1}$ is easily seen to be $n$; since (as proved by Coxeter himself) the Coxeter elements are all conjugate, they all have the same order.

Now we consider posets defined by permutations. Given $w\in\sss n$, define $\prmp w$ to be the set $\{1,\dots,n\}$, with $i\cls j$ \iff $i\ls j$ and $w(i)\ls w(j)$. Note that if $w\in\sss n$ then $\prmp w$ and $\prmp {w^{-1}}$ are isomorphic: in fact $w$ is itself an isomorphism from $\prmp w$ to $\prmp{w^{-1}}$.

Our main result in this section is the following.

\begin{thm}\label{bccox}\indent
\begin{enumerate}[ref=(\arabic{enumi})]
\vspace{-\topsep}
\item\label{bciffcox}
Suppose $w\in\sss n$. Then $\prmp w$ is a \bc \iff $w$ is a Coxeter element of $\sss n$.
\item\label{everybccox}
Every $n$-element \bc is isomorphic to $\prmp w$ for some Coxeter element $w\in\sss n$.
\item\label{isocox}
Suppose $v$ and $w$ are Coxeter elements of $\sss n$. Then $\prmp v\cong \prmp w$ \iff $v\in\{w,w^{-1}\}$.
\end{enumerate}
\end{thm}

\begin{pf}
We assume $n\gs3$. We begin with the ``only if'' part of \ref{bciffcox}. The construction of $\prmp w$ means that the number of pairs of incomparable elements in $\prmp w$ is $l(w)$. So by \cref{noofpairs} $\prmp w$ can only be a \bc if $l(w)=n-1$. If $l(w)=n-1$ but $w$ is not a Coxeter element, then $w$ can be written as a product of Coxeter generators with one Coxeter generator, say $s_m$, not occurring. This means that $w$ is contained in the \emph{Young subgroup} $\sss{(m,n-m)}$: this is the subgroup of $\sss n$ consisting of all permutations $v\in\sss n$ such that $v(\{1,\dots,m\})=\{1,\dots,m\}$. So in $\prmp w$ we have $i\cls j$ whenever $i\ls m<j$. But now the incomparability graph of $\prmp w$ is disconnected, so by \cref{bichainbiclique} $\prmp w$ is not a \bc.

So $\prmp w$ is a \bc only if $w$ is a Coxeter element. \cref{posequal,sameposets} show that every \bc is isomorphic to $\prmp w$ for some permutation $w$, which must therefore be a Coxeter element, so \ref{everybccox} is proved. Now consider the function $\psi:w\mapsto \prmp w$ from the set of Coxeter elements of $\sss n$ to the set of isomorphism classes of $n$-element posets. We have just seen that that the image of $\psi$ contains all $n$-element \bcs, and therefore this image has size at least $2^{n-3}$, by \cref{numberbcs}. On the other hand, since $w\neq w^{-1}$ but $\psi(w)=\psi(w^{-1})$ for each Coxeter element $w$, the image of $\psi$ has size at most $2^{n-3}$, by \cref{coxprops}. So equality holds everywhere; hence the image of $\psi$ is precisely the set of $n$-element \bcs (which is the ``if'' part of \ref{bciffcox}), and different inverse pairs of Coxeter elements map to different \bcs (which gives \ref{isocox}).
\end{pf}

\section{Covers and linear extensions}\label{covsec}

In this section we answer a question posed by one of the referees for this paper, showing how to count covers and linear extensions of a \bc.

First we count covers in a \bc $P$, i.e.\ pairs $a,b\in P$ such that $b$ covers $a$. Our main result is the following.

\begin{propn}\label{countcovers}
Suppose $P$ is a \bc of size $n$, and write $P=\intn{n_1}\splice\cdots\splice\intn{n_s}$ for $n_1,\dots,n_s\gs3$ with $n_1+\dots+n_s=n+2s-2$. Then the number of covers in $P$ is $2n-s-4$.
\end{propn}

\begin{pf}
We use induction on $s$: in the case $s=1$, we have $P=\intn n$, and the covers are the pairs $(a,a+2)$ and $(a,a+3)$, of which there are $2n-5$. For the inductive step it suffices to observe that for two \bcs $P$ and $Q$ the number of covers in $P\splice Q$ equals the number of covers in $P$ plus the number of covers in $Q$; this is easy to see from the construction of splices.
\end{pf}

Now we consider the number of linear extensions. For any \bc $P$, let $l(P)$ denote the number of linear extensions of $P$, i.e.\ total orders that refine the partial order on $P$. Our main result here is a recurrence relation for $l(P)$.

\begin{propn}\label{countlinear}
Suppose $P$ is a \bc with $n$ elements, and write $P=\intn{n_1}\splice\cdots\splice\intn{n_s}$ with $n_1,\dots,n_s\gs3$ and $n_1+\dots+n_s=n+2s-2$. If $n_1=\dots=n_s=3$ (so that $s=n-2$), then $l(P)=n$. Otherwise, let $r\ls s$ be maximal such that $n_r\gs4$. Then
\[
l(P)=l\bigl(\intn{n_1}\splice\cdots\splice\intn{n_{s-1}}\splice\intn{n_s-1}\bigr)+l\bigl(\intn{n_1}\splice\cdots\splice\intn{n_{r-1}}\splice\intn{n_r-2}\bigr).
\]
\end{propn}

To prove this, we need to consider the poset obtained by removing a maximal element from a \bc. Recall from \cref{twomax} that a \bc with at least three elements has exactly two maximal elements, of which exactly one is supermaximal.

We also recall some standard poset notation: we write $\ul r$ for a totally-ordered set with $r$ elements, and given posets $P$ and $Q$, we write $P\oplus Q$ for the poset obtained by taking the union of disjoint copies of $P$ and $Q$, and setting $p\cls q$ for every $p\in P$ and $q\in Q$.

\begin{lemma}\label{removemax}
Suppose $P$ is a \bc with at least three elements, and write $P=\intn{n_1}\splice\cdots\splice\intn{n_s}$ with $n_1,\dots,n_s\gs3$.
\begin{enumerate}
\item
The poset obtained by removing the supermaximal element from $P$ is isomorphic to
\[
\intn{n_1}\splice\cdots\splice\intn{n_{s-1}}\splice\intn{n_s-1}.
\]
\item
Suppose $1\ls r\ls n$ with $n_r\gs4$ and $n_{r+1}=\dots=n_s=3$, and let $Q$ be the poset obtained by removing the non-supermaximal maximal element from $P$. Then $Q$ is isomorphic to
\[
\intn{n_1}\splice\cdots\splice\intn{n_{r-1}}\splice\intn{n_r-2}\ \oplus\ \ul{s-r+1}.
\]
\end{enumerate}
\end{lemma}

\begin{pfenum}
\item
First consider the case $s=1$, so that $P=\intn{n_1}$. The supermaximal element of $\intn{n_1}$ is $n_1$, and removing this leaves the poset $\intn{n_1-1}$, as claimed. The general case follows by applying the splicing operation.
\item
First consider the case $r=s=1$, so that $P=\intn{n_1}$. The non-supermaximal maximal element of $\intn{n_1}$ is $n_1-1$; removing this leaves a poset with a greatest element $n_1$, with the other elements forming the poset $\intn{n_1-2}$. So $\intn{n_1}\setminus\{n_1-1\}\cong\intn{n_1-2}\oplus\ul1$, as claimed.

To generalise to the case $r=1\ls s$, observe that splicing $\intn{n_1}$ with $\intn3\splice\cdots\splice\intn3$ involves adding a chain of length $s-r$ lying above all elements of $\intn{n_1}$ except $n_1-1$. So the non-supermaximal maximal element is still $n_1-1$, and removing this leaves $(\intn{n_1-2}\oplus\ul1)\oplus\ul{s-r}\cong\intn{n_1-2}\oplus\ul{s-r+1}$.

Now the case $r>1$ follows by applying the splicing operation.
\end{pfenum}

\begin{pf}[Proof of \cref{countlinear}]
The case where $n_1+\dots+n_s=3$ is precisely the case where $P$ consists of a chain of size $n-1$ and an isolated element; clearly in this case $P$ has $n$ linear extensions. So assume that $n_r\gs4$ for some $r$, and let $r$ be maximal with this property.

In any linear extension of $P$, the greatest element must be one of the maximal elements of $P$. Moreover, if $p$ is a maximal element of $P$, then the number of linear extensions of $P$ having $p$ as greatest element equals the number of linear extensions of $P\setminus\{p\}$. Now the result follows from \cref{twomax,removemax}.
\end{pf}

Now we can deduce bounds for the number of covers and the number of linear extensions of a \bc. Let $F_n$ denote the $n$th Fibonacci number, with the convention that $F_0=0$. The following result follows from \cref{countcovers,countlinear} by an easy induction.

\begin{cory}\label{hilocovers}
Suppose $P$ is a \bc of size $n\gs3$, and let $c$ denote the number of covers in $P$ and $l$ the number of linear extensions of $P$. Then:
\begin{enumerate}
\item
$c\gs n-2$ and $l\gs n$, with equality in each case \iff $P\cong\intn3\splice\cdots\splice\intn3$;
\item
$c\ls2n-5$ and $l\ls F_{n+1}$, with equality in each case \iff $P\cong\intn n$.
\end{enumerate}
\end{cory}

\section{$3$-chains and beyond}\label{gensec}

In this final section we briefly discuss a natural generalisation of \bcs. Given $r\in\bbn$, define an \emph{$r$-chain} to be a poset $P$ such that $P$ is uniquely expressible as the union of $r$ chains, and the partial order on $P$ is maximal with this property.

For a first family of examples, we can generalise the posets $\intn n$ defined in \cref{bcsec}. Given $n\in\bbn$, define $\intn n^r=\{1,\dots,n\}$, with $i\cls j$ \iff $i\ls j-r$ in the usual order on $\{1,\dots,n\}$.

To construct all $r$-chains, it seems to be possible to generalise the recursive construction suggested after \cref{autotriv}: starting from the unique $r$-element $r$-chain, repeatedly add supermaximal elements, at each stage choosing which $r-1$ of the existing $r$ maximal elements should remain maximal. From this it should follow that the number of isomorphism classes of \emph{labelled} $r$-chains $P$ with $n\gs r$ elements (i.e.\ with the $r$ chains labelled $1,\dots,r$, and with isomorphisms required to preserve labelling) is $r^{n-r}$.

However, the other results in this paper remain to be extended. The author cannot see an obvious way to extend the constructions in \cref{consec} to the cases $r\gs3$. For the results in \cref{splicesec}, it should be possible to splice $r$-chains to make larger ones, but a direct analogue of \cref{splicepn} fails: there are $r$-chains other than $\intn1^r,\intn2^r,\dots$ which cannot be written as splices of smaller $r$-chains. The first example with $r=3$ is as follows.
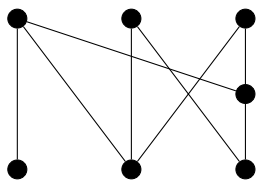
\[
\begin{tikzpicture}[scale=1]
\nd{0,0}a
\nd{0,2}b
\nd{1.5,0}c
\nd{1.5,2}d
\nd{3,0}e
\nd{3,1}f
\nd{3,2}g
\draw(a)--(b)--(c)--(d)--(e)--(f)--(g)--(c);
\draw(b)--(f);
\end{tikzpicture}
\]
The results of \cref{coxsec} also have no obvious extension beyond $r=2$; indeed, $r$-chains need not have dimension $2$ when $r>2$, so not all $r$-chains have the form $\prmp w$ for $w$ a permutation. For example, $\intn7^3$ has dimension $3$.

We hope to say more about $r$-chains in future work.

\section{Infinite \bcs}\label{infinitesec}

In this section we briefly consider infinite \bcs. The definition of a \bc still makes sense for infinite posets, and many of the results in this paper appear to have analogues, provided we restrict attention to \emph{locally finite} posets, i.e.\ posets $P$ such that for any two elements $p,q$, there are only finitely many elements $r$ with $p\cls q\cls r$. So for this section we assume all posets are locally finite, and we explain how to generalise some of the results in this paper. Some of the statements made here are conjectural; we expect they should not be hard to prove, though of course the inductive approach used in this paper for finite \bcs will not work.

The simple examples of \bcs in \cref{bcsec} generalise naturally: we can take any locally finite chain (which may have a greatest or a least element, or neither) together with an isolated element to give an infinite \bc. We can also construct analogues of the finite \bcs $\intn n$, by taking the underlying set to be $\bbn$, $-\bbn$, or $\bbz$, with $i\cl j$ when $i\ls j-2$. We call the resulting \bcs $\intn{+\infty}$, $\intn{-\infty}$, $\intn{\infty}$ respectively.

Now we consider maximal elements. An infinite \bc has at most two maximal elements; it has a supermaximal element (which is unique) \iff it has two maximal elements, and this in turn happens \iff every element lies below a maximal element. A corresponding statement holds for minimal elements. As a consequence, an infinite \bc cannot have two maximal and two minimal elements. We refer to a \bc with at most one maximal element and at most one minimal element as a \emph{doubly-infinite} \bc.

The construction using binary sequences in \cref{consec} can be generalised to the infinite case, though with slight complications. First we consider infinite \bcs with two minimal elements. To give a construction of these, take an infinite binary sequence $a=a_1a_2\dots$. Define an infinite family of sequences $a(0),a(1),\dots$ by
\[
a(n)=a_1\,\dots\,a_{n-1}\,0\,1\,a_{n+1},a_{n+2},\,\dots
\]
(with the $0$ omitted in the case $n=0$). Additionally, if $a_r=0$ for all $r\gg0$, define the sequence $a(\infty)=a$. The sequences defined in this way are all the sequences that can be obtained by inserting a $0$ or a $1$ into $a$, and the poset $\calp_a$ can be defined as in \cref{consec}, and is a \bc with two minimal elements (namely, the sequences obtained by adding a $0$ or a $1$ at the start of $a$). $\calp_a$ has a maximal element (namely, $a$ itself) \iff $a_r$ is constant for $r\gg0$.

A dual construction (using backwards-infinite sequences $\dots a_{-2}a_{-1}$) works for infinite \bcs with two maximal elements. For doubly-infinite \bcs, we need to use doubly-infinite binary sequences $a=\dots a_{-1}a_0a_1a_2\dots$. Given such a sequence, we define sequences $a(n)$ for $n\in\bbz$ by
\[
a(n)_r=
\begin{cases}
a_r&(r<n)\\
0&(r=n)\\
1&(r=n+1)\\
a_{r-1}&(r>n+1).
\end{cases}
\]
In addition:
\begin{itemize}
\item
if $a_r=0$ for all $r\gg0$, define the sequence $a(\infty)=a$;
\item
if $a_r=1$ for all $r\ll0$, define the sequence $a(-\infty)$ to be the sequence $\shift a$ obtained from $a$ by shifting entries by $1$, i.e.\ $\shift a_r=a_{r-1}$ for all $r$.
\end{itemize}
The sequences constructed are then all the sequences that can be obtained from $a$ by inserting a symbol in $a$ and shifting by $1$ every symbol to the right of the inserted symbol. Now we can follow the construction of \cref{consec}, yielding an infinite \bc $\calp_a$. This \bc has a maximal element (namely $a$) \iff $a_r$ is constant for $r\gg0$, and has a minimal element (namely $\shift a$) \iff $a_r$ is constant for $a\ll0$. It seems likely that these constructions give all locally finite \bcs.

Locally finite \bcs can also be constructed by splicing, using the finite \bcs $\intn n$ as well as the \bcs $\intn{+\infty}$, $\intn{-\infty}$ and $\intn{\infty}$ introduced above. We suspect that every infinite \bc is obtained by splicing \bcs $\intn n$. More specifically, every infinite \bc with two minimal elements is either of the form
\[
\intn{n_1}\splice\intn{n_2}\splice\cdots
\]
for some infinite sequence $n_1,n_2,\dots$ of integers greater than $2$, or of the form
\[
\intn{n_1}\splice\cdots\splice\intn{n_s}\splice\intn{+\infty}
\]
for some finite sequence $n_1,\dots,n_s$. Similarly every doubly-infinite \bc should have one of the forms
\[
\cdots\splice\intn{n_{-1}}\splice\intn{n_0}\splice\intn{n_1}\splice\cdots,
\]
\[
\cdots\splice\intn{n_{-2}}\splice\intn{n_{-1}}\splice\intn{+\infty},
\]
\[
\intn{-\infty}\splice\intn{n_1}\splice\intn{n_2}\splice\cdots,
\]
\[
\intn{\infty}
\]
for an appropriate (finite or infinite) sequence of integers $n_r\gs3$. One way to prove this would be as follows: define a \emph{super} element in a \bc $P$ to be an element $p$ which is comparable with every element of $P$ except one (say $p^\ast$). Then the ideal generated by $p$ and $p^\ast$ should be a \bc, and the coideal generated by $p$ and $p^\ast$ (that is, the set $\lset{a\in P}{p\cls a\text{ or }p^\ast\cls a}$) should be a \bc, and $P$ is the splice of these two \bcs. The main task then is to prove that a \bc with no super elements is isomorphic to $\intn\infty$ (and similarly that a \bc whose only super element is minimal is isomorphic to $\intn{+\infty}$).

We expect the results on graphs in \cref{graphsec} to generalise to infinite \bcs; the incomparability graphs of \bcs should then be caterpillars in which all vertices have finite degree except possibly the head and the tail (i.e.\ the endpoints (if there are any) of the path formed by the non-leaves).

Finally, we expect that there should be analogues of the results in \cref{coxsec}, for suitable infinite analogues of Coxeter elements. We begin with the case of \bcs having two minimal elements. Let $\sss{\bbn}$ denote the group of permutations of $\bbn$, and define a \emph{Coxeter element} in $\sss{\bbn}$ to be a permutation $w$ with the property that for every $r>1$, one of $w(r)$, $w^{-1}(r)$ is greater than $r$ and the other is less than $r$. Now the poset $\prmp w$ defined as in \cref{coxsec} should be a \bc with two minimal elements and no maximal elements. For \bcs with two minimal elements and a maximal element, a more complicated notion is needed: rather than permutations, one should take bijections from $\bbn$ to $\bbn\cup\{\infty\}$. For doubly-infinite \bcs, we make a similar definition: define a Coxeter element in $\sss{\bbz}$ to be a permutation of $\bbz$ such that for every $r$ one of $w(r)$ and $w^{-1}(r)$ is greater than $r$ and the other is less than $r$. Then the poset $\prmp w$ should be a \bc with no maximal or minimal elements. To obtain doubly-infinite \bcs with a maximal or minimal element, one needs to consider elements $\pm\infty$; we leave the reader to work out the details.


\begin{thebibliography}{99}

\backrefparscanfalse

\bibi{0hecke}{F}
{M.~Fayers}
{$0$-Hecke algebras of finite Coxeter groups}
{J.\ Pure Appl.\ Algebra}
{199}{2005}{27--41}

\bibi{gallai}{G}
{T.~Gallai}
{Transitiv orientierbare Graphen}
{Acta Math.\ Acad.\ Sci.\ Hungar.}
{18}{1967}{25--66}

\bibi{hasc}{HS}
{F.\ Harary \& A.\ Schwenk}
{The number of caterpillars}
{Disc.\ Math.}6{1973}{359--365}

\bibbook{hum}{H}
{J.\ Humphreys}
{Reflection groups and Coxeter groups}
{Cambridge Studies in Advanced Mathematics}{29}{Cambridge University Press, Cambridge}{1990}

\bibi{prt}{P}
{O.~Pretzel}
{A representation theorem for partial orders}
{J.\ London Math.\ Soc.}{42}{1967}{507--508}

\bibitem[R]{mo}
J.~Rickard, answer to `Has anyone seen these posets before?', mathoverflow.net/q/278244.\backrefprint\renewcommand\con{}\renewcommand\cons{}

\bibi{shi}{S}
{J.--Y.~Shi}
{The enumeration of Coxeter elements}
{J.~Algebraic Combin.}
6
{1997}
{161--171}

\end{thebibliography}
\end{document}